\documentclass[12pt,a4paper]{article}

\IfFileExists{ajr.sty}{
\usepackage[mathbf,mathcal]{euler}\usepackage{beton}
\usepackage{defns}
\usepackage[hyperref,colour]{ajr}
}{ \usepackage{defns} }

\newcommand{\pde}{{\textsc{pde}}}
\newcommand{\jp}{_{j+1}}
\newcommand{\jm}{_{j-1}}
\newcommand{\jpm}{_{j\pm1}}
\newcommand{\jmm}{_{j-2}}
\newcommand{\jpp}{_{j+2}}

\newcommand{\cM}{{\mathcal M}}
\newcommand{\cA}{{\mathcal A}}
\newcommand{\cB}{{\mathcal B}}

\makeatletter
\def\arcsinh{\mathop{\operator@font arcsinh}\nolimits}
\makeatother
\newcommand{\ibc}{\textsc{ibc}}

\begin{document}

\title{A holistic finite difference approach models linear dynamics
consistently}
\author{AJ Roberts\thanks{Dept.\ Maths \& Comput, University of
Southern Queensland, Toowoomba, Queensland 4352, \textsc{Australia}.
\protect\url{mailto:aroberts@usq.edu.au}}}
\maketitle

\begin{abstract}
    I prove that a centre manifold approach to creating finite
    difference models will consistently model linear dynamics as the
    grid spacing becomes small.  Using such tools of dynamical systems
    theory gives new assurances about the quality of finite difference
    models under nonlinear and other perturbations on grids with
    finite spacing.  For example, the advection-diffusion equation is
    found to be stably modelled for all advection speeds and all grid
    spacing.  The theorems establish an extremely good form for the
    artificial internal boundary conditions that need to be introduced
    to apply centre manifold theory.  When numerically solving
    nonlinear partial differential equations, this approach can be
    used to derive systematically finite difference models which
    automatically have excellent characteristics.  Their good
    performance for finite grid spacing implies that fewer grid points
    may be used and consequently there will be less difficulties with
    stiff rapidly decaying modes in continuum problems.
\end{abstract}

\paragraph{Maths.~Subj.~Class:} 37L65, 65M20, 37L10, 65P40, 37M99

\tableofcontents

\section{Introduction}

Following the introduction of holistic finite differences
in~\cite{Roberts98a, MacKenzie00a}, we would like to investigate
numerical models for the dynamics of a field~$u(x,t)$ evolving
according to a nonlinear reaction-diffusion equation such as
$u_t=u_{xx}+f(u,u_x)$\,.  This particular class includes Burgers'
equation, $f=-uu_x$, and autocatalytic reactions, such as $f=u(1-u)$.
However, before attacking such nonlinear problems,  here we restrict
attention to proving that the new methodology accurately models the
dynamics of quite general \emph{linear} \pde's.

Modern dynamical systems theory has had to date very little impact on
classical numerical approximations.
Indeed, the very first sentence in Garc\'ia-Archilla
\etal~\cite{Archilla98} says ``Finite-element methods seem not to have
benefited as much as spectral methods from some of the recent advances
in the Dynamical Systems approach to partial differential equations''.
The concept of inertial manifolds has been developed to capture the
long-term, low-dimensional dynamics of dissipative \pde's
\cite{Temam90}.
However, most efforts to construct approximations to an inertial
manifold have been based upon the \emph{global} nonlinear Galerkin
method of Roberts~\cite{Roberts89}, Foias \etal~\cite{Foias88c} and
Marion \& Temam \cite{Marion89}.
This is so even for the variants explored by Jolly
\etal~\cite{Jolly90} and Foias \cite{Foias91b}.
In contrast, the approach proposed here is based purely upon the
\emph{local} dynamics on small elements while maintaining, as do
inertial manifolds, fidelity with the solutions of the original \pde.

\begin{figure}[tbp]
\begin{center}
\setlength{\unitlength}{1ex}
\begin{picture}(50,48)
    \put(15,15){
    \put(0,0){\vector(1,0){30}}\put(30.5,-.5){$h$}
    \put(0,0){\vector(0,1){30}}\put(-6,30.5){nonlinearity}
    \put(0,0){\vector(-1,-2){6.5}}\put(-8,-14.5){$\gamma$}
    \put(-2,-.5){$0$}
    \put(0,25){\line(1,0){20}}
    \put(20,0){\line(0,1){25}}
    {\thicklines \put(20,0){\circle{2}}}
    \put(-5,-10){
    \put(-8,-.5){$\gamma=1$}
    \put(0,0){\line(0,1){25}}
    \put(20,0){\line(0,1){25}}
    \put(20,25){\circle*{2}}
    \put(21.5,25){physical problem}
    \put(20,0){\circle*{2}}
    \put(22,-2){\parbox{14ex}{\raggedright general linear problem}}
    \put(4,0.5){\it consistency}
    \put(4,25.5){\it consistency}
    \thicklines
    \put(0,0){\circle{2}}
    \put(0,0){\vector(1,0){19}}
    \put(0,25){\circle{2}}
    \put(0,25){\vector(1,0){19}}
    } 
    \put(0,25){\line(-1,-2){5}}
    \put(20,25){\line(-1,-2){5}}
    \put(18.5,5){\it holistic}
    \put(18.5,-5){\it holistic}
    \thicklines
    \put(20,0){\vector(-1,-2){4.5}}
    \put(20,0){\vector(-1,3){4.6}}
    }
\end{picture}
\end{center}
    \caption{conceptual diagram showing: the traditional finite
    difference modelling approaches (rightward-arrows) the physical
    problem (upper disc) via asymptotic consistency as the grid size
    $h\to 0$ (left circles); whereas the holistic method approaches
    (forward-arrows) the physical problem via asymptotics in
    nonlinearity and the inter-element coupling $\gamma$ (from right
    circle).  Herein we establish how to use the holistic approach to
    do \emph{both} in order to model a general linear problem (lower
    disc).}
    \label{fig:conc}
\end{figure}
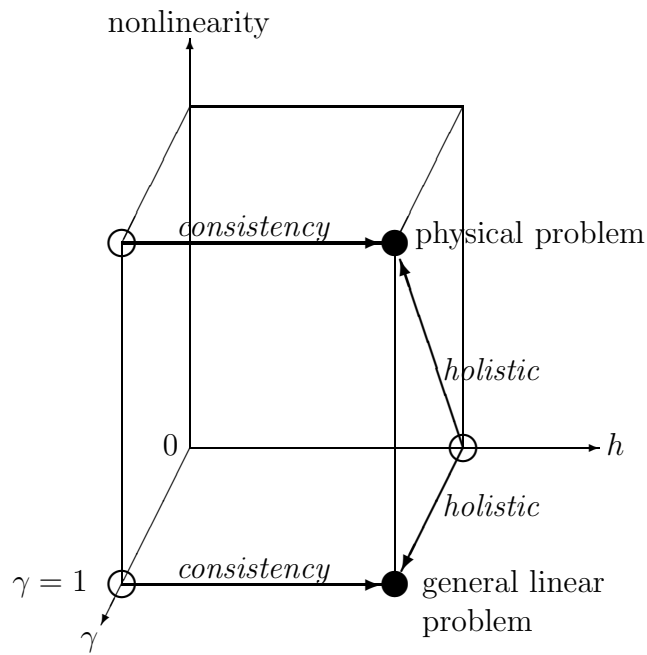

I propose \cite{Roberts98a} to use centre manifold theory to construct
finite difference models.
For problems in one spatial dimension consider implementing the method
of lines by discretising in $x$ and integrating in time as a set of
ordinary differential equations, sometimes called a semi-discrete
approximation \cite[e.g.]{Foias91,Foias91b}.
We only address spatial discretisation and treat the resulting set of
ordinary differential equations as a continuous time dynamical system.
Classical finite difference approximations are made by appealing to
consistency in the limit as the grid spacing $h\to0$; traditionally
one constructs models to errors $\Ord{h^2}$ or $\Ord{h^4}$ depending
upon small $h$ asymptotics, shown schematically by the
rightward-arrows in Figure~\ref{fig:conc}.
In contrast, we here analyse the dynamics at fixed grid spacing $h$
and use centre manifold theory to accurately model the nonlinear
dynamics---theory~\cite[e.g.]{Carr81} assures us that the
low-dimensional, numerical model then accurately captures the dynamics
in an expansion in the nonlinearity, shown schematically as the
forward-arrows in Figure~\ref{fig:conc}.
The analysis rests upon the exponential decay of the small, subgrid
structures in each local element.
Being essentially local in space, the analysis here is flexible enough
to subsequently cater for spatial boundaries and spatially varying
coefficients.
I call the model ``holistic'' because the centre manifold is made up
of actual \emph{solutions} of the \pde{} and is thus invariant under
algebraic rewriting of the governing equations.
However, to apply the centre manifold theory we have to use a homotopy
in a parameter~$\gamma$: when $\gamma=0$ the discrete finite elements
of the domain are completely uncoupled from each other; when
$\gamma=1$, the requisite continuity is reclaimed and the physical
\pde{} solved.
The caveat is that the centre manifold model has to be used at
$\gamma=1$ whereas the supporting theory only guarantees accuracy in a
neighbourhood of $\gamma=0$; we aim to make the useful neighbourhood
big enough to include the relevant $\gamma=1$ (this sort of technique
has proven effective in thin fluid flows \cite[e.g.]{Roberts94c}).
One way to reasonably secure the centre manifold model, and the way
explored herein, is to require that the model is \emph{also
consistent} with the \pde{} as the grid spacing $h\to 0$.
Thus we aim to construct finite difference numerical models that not
only are justified by their asymptotic expansions in nonlinearity and
$\gamma$, but are also justified by an asymptotic expansion in $h$
(see Figure~\ref{fig:conc}).
This dual justification is the completely novel feature of the
approach.

The first step, taken in this paper, is to establish a centre manifold
approach that is also guaranteed to construct a consistent finite
difference model for a general \emph{linear} \pde, shown schematically
in Figure~\ref{fig:conc} as the disc in the $\gamma h$-plane (zero
nonlinearity).
We leave to later research the problem of guaranteeing the consistent
modelling of nonlinear dynamics.
Herein we explore the finite difference modelling of the linear \pde{}
\begin{equation}
    \D tu=\cA u+\epsilon \cB u\,,
    \label{eq:pde}
\end{equation}
where the linear operator $\cA$, presumed generally dissipative, is
even (it contains only even order derivatives in $x$ with constant
coefficients), and $\cB$ is an odd linear operator (it contains only
odd order derivatives with constant coefficients); $\cA$ is assumed
dissipative for all modes except $u=\mbox{const}$.
The case of space-time varying coefficients to the linear problem is
also left for later study; however, expect that because the analysis
here is local in~$x$, then such varying coefficients can be treated as
a perturbing influence to the basic analysis herein.
As a specific example, we discuss in \S\ref{Sadvec} the linear
advection-diffusion equation $u_t=-\epsilon u_x+u_{xx}$ and discover
many remarkable properties of the holistic finite difference models.
The separation between the two types of linear terms into $\cA u$ and
$\epsilon\cB u$ occurs because first we prove consistency for even
terms, \S\ref{Sline}, before moving on to prove consistency for the
odd terms, \S\ref{Slino}.
The results are verified for the perturbed diffusion equation by the
computer algebra program listed in the Appendix.

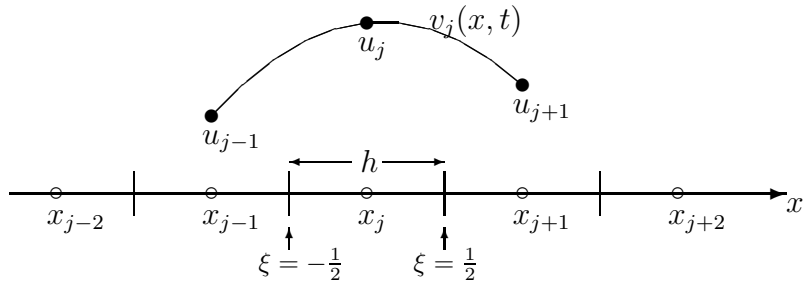
\begin{figure}
\begin{center}
    \setlength{\unitlength}{0.1em}
    \begin{picture}(260,82)
          \thinlines
	  \put(65,50){
	  \put(0,0){\line(1,1){10}}
	  \put(10,10){\line(5,4){10}}
	  \put(20,18){\line(5,3){10}}
	  \put(30,24){\line(5,2){10}}
	  \put(40,28){\line(5,1){10}}
	  \put(50,30){\line(5,0){10}}
	  \put(60,30){\line(5,-1){10}}
	  \put(70,28){\line(5,-2){10}}
	  \put(80,24){\line(5,-3){10}}
	  \put(90,18){\line(5,-4){10}}
	  \put(0,0){\circle*{4}} \put(-3,-9){$u\jm$}
	  \put(50,30){\circle*{4}} \put(47,21){$u_j $}
	  \put(100,10){\circle*{4}} \put(97,1){$u\jp$}
	  \put(70,28){$v_j(x,t)$}
	  }
	  \put(90,7){\vector(0,1){7}}
	  \put(140,7){\vector(0,1){7}}
	  \put(80,0){\footnotesize$\xi=-\half$}
	  \put(130,0){\footnotesize$\xi=\half$}
	  \put(113,32){$h$}
	      \put(120,35){\vector(1,0){20}}
	      \put(110,35){\vector(-1,0){20}}
          \put(12,15){$x\jmm$}
          \put(62,15){$x\jm$}
          \put(112,15){$x_j$}
          \put(162,15){$x\jp$}
          \put(212,15){$x\jpp$}
	  \multiput(40,18)(50,0){4}{\line(0,1){14}}
	  \multiput(15,25)(50,0){5}{\circle{3}}
	  \thicklines
	  \put(250,19){$x$}
	  \put(0,25){\vector(1,0){250}}
    \end{picture}
\end{center}
\caption{schematic picture of the equi-spaced grid, $x_j$ spacing $h$,
the unknowns~$u_j$, the artificial internal boundaries between each
element (vertical lines), and in the neighbourhood of $x_j$ the field
$v_j(x,t)$ which extends outside the element and, if $\gamma=1$,
passes through the neighbouring grid values $u\jpm$.}
\label{fig:grid}
\end{figure}

Introduce a regular grid as shown in Figure~\ref{fig:grid} with grid
points a distance $h$ apart, $x_j=jh$ for example, and using $u_j $ to
denote the value of the field at each grid point,
\begin{equation}
    u_j(t)=u(x_j,t)\,.
    \label{eq:amp}
\end{equation}
We express the field in the neighbourhood of the $j$th grid point by
$u=v_j(x,t)$.
We do not restrict the function $v_j$ to just the $j$th element, but
allow it to extend analytically out to at least the adjacent grid
points as shown in Figure~\ref{fig:grid}.

Herein we establish the small $h$ consistency that follows from using
the nonlocal, internal ``boundary conditions''
\begin{equation}
    v_j(x\jpm,t)=(1-\gamma)v_j(x_j,t)+\gamma v\jpm(x\jpm,t)\,.
    \label{eq:ibcv}
\end{equation}
That is, the field of the $j$th element when evaluated at the
surrounding gridpoints, $v_j(x\jpm,t)$, is a continuation between two
critical extremes: when $\gamma=1$, it is the field at those grid
points, $v\jpm(x\jpm,t)$, to in effect recover the physical continuity
as shown schematically in Figure~\ref{fig:grid}; but when $\gamma=0$,
the field is just identical to the mid-element value $v_j(x_j,t)$ so
that the element becomes isolated from all neighbours.
Equivalently, (\ref{eq:ibcv}) is transformed to the following
appealing two difference equations\footnote{Natural symmetry also
makes the general results most easy to express in terms of centred
difference operators and so I use them throughout this paper.}
\emph{evaluated at the centre of the element}, $x=x_j$,
\begin{equation}
    \mu_x\delta_x v_j(x,t)=\gamma\mu\delta v_j(x_j,t)
    \quad\mbox{and}\quad
    \delta_x^2 v_j(x,t)=\gamma\delta^2 v_j(x_j,t)\,.
    \label{eq:jbcv}
\end{equation}
That is, in the two extremes: when $\gamma=0$ the first and second
differences have to be zero; whereas when $\gamma=1$ the first two
differences of the field centred on each element have to agree
with the first two differences of the grid values.
Note a distinction which is very important throughout this work:
unadorned difference operators, such as the central mean
$\mu=(E^{1/2}+E^{-1/2})/2$ and central difference
$\delta=E^{1/2}-E^{-1/2}$ written in terms of the shift operator
$Eu_j=u_{j+1}$ \cite[p64,e.g.]{npl61}, apply to the grid index~$j$
(with step~$1$) whereas those with subscript $x$, as in $\mu_x$ and
$\delta_x$, are differences in $x$ only (with step~$h$).
Using the definition of the amplitudes~(\ref{eq:amp}), these internal
boundary conditions (\ibc) simplify to the following form which we use
throughout \S\S\ref{Sadvec}--\ref{Slino}: evaluated at $x=x_j$
\begin{equation}
    \mu_x\delta_x v_j(x,t)=\gamma\mu\delta u_j
    \quad\mbox{and}\quad
    \delta_x^2 v_j(x,t)=\gamma\delta^2 u_j\,.
    \label{eq:jbc}
\end{equation}
In actually developing finite difference models these \ibc's may take
any of many equivalent forms \cite[e.g.]{Roberts98a}.
Small~$h$ consistency seems easiest to establish in this particular
discrete form.

\section{The dynamics collapses onto a centre manifold}
\label{Scm}

I establish here the basis of a centre manifold analysis of the linear
\pde~(\ref{eq:pde}).
It appears necessary, and is the route taken here, to separate the
linear effects into those generated by even terms, represented by
$\cA$ and presumed generally dissipative on the grid scale~$h$ but
with one $0$~eigenvalue corresponding to
$u=\mbox{const}$,\footnote{The analysis presented here could be
applied to modelling unstable dynamics, such as that from negative
diffusion $u_t=-u_{xx}$.
The difference is that the relevance theorem would no longer
apply---$\cM$ would not be attractive and the finite difference model
would not capture the long term dynamics.
The centre manifold model would, however, capture all the finite
solutions, if any.} and those generated by odd terms, represented by
$\epsilon\cB$.
The analysis is to be based on the situation when the coefficient of
the odd terms $\epsilon=0$ and when adjacent elements are decoupled,
$\gamma=0$ (the right-hand circle in Figure~\ref{fig:conc}).
Then centre manifold theory \cite[e.g.]{Carr81} guarantees the
existence and relevance of a numerical model parametrised by the
discrete values of the field at the grid points, $u_j$.
The constructed model is accurate to the order of the residuals of the
differential equation~(\ref{eq:pde}).

The spectrum of the linear dynamics is used to show there exists a
centre manifold.
Adjoin to~(\ref{eq:pde}) the trivial dynamical equations
\begin{equation}
    \dot\gamma=\dot\epsilon=0\,,
    \label{eq:triv}
\end{equation}
so that terms multiplied by $\epsilon$ or $\gamma$ become ``nonlinear
terms'' in the asymptotic expansion we develop about $\gamma =\epsilon
=0$.
Setting $\epsilon=0$ eliminates the odd terms to leave simply the
linear equation $u_t=\cA u$; for example, the diffusion equation
$u_t=u_{xx}$.
This is to be solved in the vicinity of $x_j$ for a field
$u=e^{\lambda t}w(x)$ where here the dependence upon $j$ is implicit
in the eigenfunction~$w$ of $\cA w=\lambda w$.
This is a constant coefficient \ode{} and so has trigonometric general
solutions $w\propto e^{i\alpha x}$ with corresponding eigenvalue
$\lambda(\alpha)$ which is negative for non-zero wavenumber $\alpha$
as $\cA$ is presumed generally dissipative; for example,
$\lambda=-\alpha^2$ for the diffusion equation.
The appropriate boundary conditions come from the nonlocal decoupling
conditions that $w(x\jpm)=w(x_j)$ from~(\ref{eq:ibcv}) with $\gamma=0$.
Thus within each element the eigenmodes are: $\exp[i\alpha_n (x-x_j)]$
for even integer $n$, where $\alpha_n=n\pi/h$; and also
$\sin[\alpha_n(x-x_j)]$ for odd integer $n$.
The spectrum is then $\lambda_n=\lambda(n\pi/h)$; for example,
$\lambda_n=-n^2\pi^2/h^2$ for the diffusion equation.
Thus linearly and in the absence of inter-element coupling, generally
expect all modes to decay exponentially quickly to zero on a time
scale $\Ord{h^2}$ as $\lambda(\alpha)$ will be symmetric, except for
the neutral mode, $n=0$, which is constant in~$x$, $v_j(x,t)=u_j$.
Thus for small enough $\epsilon$ and $\gamma$, theory
(\cite[p281]{Carr83b} or~\cite[p96]{Vanderbauwhede89}) assures us that
there exists a centre manifold~$\cM$ for the system~(\ref{eq:pde})
coupled across elements by~(\ref{eq:ibcv}).
The centre manifold~$\cM$ is here, by~(\ref{eq:amp}), to be
parametrised by the values of the field at the grid points, $u_j$.
Thus using $\vec u$ to denote the set of grid values $u_j$, the
``amplitudes'', theory \cite{Carr81, Carr83b, Vanderbauwhede89}
supports our description of the centre manifold and the evolution
thereon as
\begin{equation}
    u_j=v(\vec u,x)\,,\quad\mbox{such that}\quad
    \dot u_j=g(\vec u)\,,
    \label{eq:mod}
\end{equation}
where the fact that the right-hand side functions depend upon the
element~$j$ is implicit (no confusion need arise because the
translational invariance in $x$ leads to identical expressions in each
element except for appropriate changes of the subscript~$j$).
The evolution $\dot u_j=g(\vec u)$ forms the holistic finite
difference model.
When the model is constructed to errors $\Ord{\gamma^\ell}$, then we
account for interactions among $\ell-1$ elements on either side of any
given element and so the resulting finite difference model has a
stencil of width $2\ell-1$ on the spatial grid.
I call these models ``holistic'' because, unlike traditional finite
difference modelling which just analyses separately each term in the
equations, here $\cM$ is made up of actual solutions of the \pde{} and
so here the discretisation models all the possible interactions
between all the terms in the equations \cite[e.g.]{Roberts98a,
MacKenzie00a}.

Moreover, the evolution on the centre manifold forms an accurate
low-dimensional model of the dynamics of the \pde~(\ref{eq:pde}).
Again provided $\epsilon$ and $\gamma$ are small enough, theory,
\cite[p282]{Carr83b} or~\cite[p128]{Vanderbauwhede89}, assures us that
all solutions sufficiently near the centre manifold~$\cM$ are not only
attracted to $\cM$ but exponentially quickly approach the actual
solutions of the \pde{} that make up $\cM$.
The rate of attraction is approximated for practical purposes by the
leading negative eigenvalue; for example, $\lambda_1=-\pi^2/h^2$ for
the diffusion equation.
In the development of inertial manifolds by Temam~\cite{Temam90} and
others, this property is sometimes termed the asymptotic completeness
of the model, for example see Robinson~\cite{Robinson96} or Constantin
\etal~\cite[\S12-3]{Constantin89}, and sometimes as exponential
tracking \cite[e.g.]{Foias89}.
Observe that this is one of the crucial new aspects brought to finite
difference modelling by centre manifold theory.
It asserts that on a \emph{finite} grid spacing we will faithfully
track the solutions of the original \pde.
There will be some limitations: steep gradients and large
nonlinearities will test the model as always.
But the centre manifold theory, as seen here in \S\ref{Sadvec} and in
introductory work \cite{Roberts98a, MacKenzie00a}, provides a
rationale and a method to construct the requisite adjustments to a
finite difference model to robustly model a wide range of dynamics on
a finite grid spacing.
The main limitation is the rate of attraction to~$\cM$: the centre
manifold model should be able to capture any dynamics occuring on a
time-scale larger than $1/|\lambda_1|$ ($h^2/\pi^2$ for diffusion).
Thus the model evolution on $\cM$, (\ref{eq:mod}), captures all the
long term dynamics with some provisos.

\section{Advection-diffusion is modelled robustly}
\label{Sadvec}

Here we explore perhaps the simplest nontrivial example in the class
of \pde's: the advection-diffusion equation
\begin{equation}
    u_t=-\epsilon u_x+u_{xx}\,,
    \label{eq:advec}
\end{equation}
where $\epsilon$ is the advection speed that will be treated as small
in the asymptotics but investigated over a range of sizes.
This \pde{} fits into the scheme outlined in the previous sections
with the operators $\cA=\partial_x^2$ and $\cB=-\partial_x$.
We show consistency for small grid spacing $h$, and find interesting
and stable upwind approximations for large $\epsilon h$.
The associated sophisticated dependence upon $\epsilon$ is perhaps
indicative of the need to treat odd operators differently to even.

The centre manifold models discussed here were all derived by the
computer algebra program given in Appendix~\ref{Scompalg}.
The listing in the appendix replaces the recording of tedious
intermediate algebraic steps.
Suffice to say that the algorithm is based upon an iterative
refinement technique reported in~\cite{Roberts96a}: given an
approximation to the centre manifold and the evolution
thereon~(\ref{eq:mod}), the computer algebra calculates a correction
to reduce the residuals of the governing equations, the
\pde~(\ref{eq:advec}), the \ibc's~(\ref{eq:jbc}) and the amplitude
condition~(\ref{eq:amp}).
Iteration then drives the residuals to zero to some order of error.
Thus by the Approximation theorem~\cite[e.g.]{Carr81,Zhenquan00} we
are sure that the model is correct to the same order of accuracy.
In this section we proceed to use the results without any further
description of the algebra involved.

First explore the holistic finite difference model with only
first-order interactions between adjacent elements, that is, the case
$\ell=2$ where errors are $\Ord{\gamma^2}$.
As mentioned earlier, we find odd derivatives in $x$ cause a hierarchy
of refinements parametrised by increasing powers of $\epsilon$.
To low-order in $\epsilon$ the computer algebra shows the evolution on
the centre manifold is
\begin{equation}
    \dot u_j=\gamma\left[ -\frac{\epsilon}{h}\mu\delta
    +\frac{1}{h^2}\delta^2 +\frac{\epsilon^2}{12}\delta^2 \right]u_j
    +\Ord{\epsilon^3,\gamma^2}\,;
    \label{eq:lowad}
\end{equation}
substitute $\gamma=1$ to obtain a model for the advection-diffusion
\pde~(\ref{eq:advec}).
Introduced automatically in this analysis is the novel term involving
the square of the advection speed $\frac{\epsilon^2}{12}\delta^2 u_j$.
One alternative is to view this term as an upwind correction to the
finite difference approximation of the first derivative:
\begin{displaymath}
    -\frac{\epsilon}{h}\mu\delta +\frac{\epsilon^2}{12}\delta^2
    = - \frac{\epsilon}{h} \left[ \left( \frac{1}2 -\frac{\epsilon
    h}{12}\right)E^{1/2} + \left( \frac{1}2 +\frac{\epsilon
    h}{12}\right)E^{-1/2}
    \right]\delta\,,
\end{displaymath}
which increases the weight of the upwind grid point ($E^{-1/2}$ is
upwind if $\epsilon$ is positive).
Such upwind corrections are well known to be stabilising for finite
advection speeds~$\epsilon$.
The other alternative is to view the novel term as increasing the
dissipation operator, to
\begin{displaymath}
    \frac{1}{h^2}\left[ 1+\frac{\epsilon^2 h^2}{12} \right]\delta^2 \,,
\end{displaymath}
and thus also improving the stability of the finite difference scheme
for finite speeds~$\epsilon$.
It is this latter view that seems easiest to establish at higher order.

\begin{figure}[tbp]
    \centering
    \includegraphics{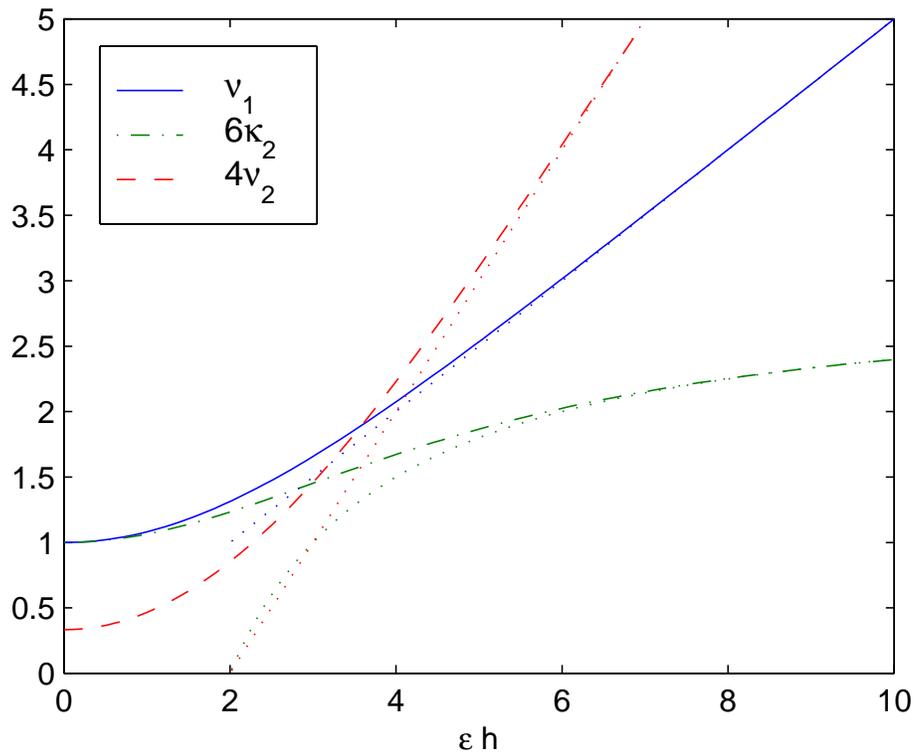}
    \caption{coefficients of the centre manifold
    models~(\ref{eq:hiad}) and~(\ref{eq:hiiad}) as a function of
    advection speed and grid spacing, $\epsilon h$.  These curves are
    at least of graphical accuracy and are obtained via the Shanks
    transform of the Taylor series to errors $\Ord{\epsilon^{14}}$.
    The dotted lines are the presumed large $\epsilon h$ asymptotes:
    $\nu_1\approx \frac{\epsilon h}{2}$\,; $\nu_2\approx
    \frac{\epsilon h}{4}-\frac{1}{2}$\,; $\kappa_2\approx
    \frac{1}{2}-\frac{1}{\epsilon h}$\,.}
    \label{fig:egcof}
\end{figure}

It is interesting to explore in some detail to high order in
$\epsilon$.
I find the model is (after setting $\gamma=1$)
\begin{eqnarray}
    \dot u_j&=&
-{\epsilon}\frac{\mu\delta}{h} u_j
+\nu_1\frac{\delta^2}{h^2} u_j\,,
\label{eq:hiad}\\
\mbox{where}\quad\nu_1 &=& 1
+\frac{\epsilon^2h^2}{12}
-\frac{\epsilon^4h^4}{720}
+\frac{\epsilon^6h^6}{30240}
-\frac{\epsilon^8h^8}{1209600} +\Ord{\epsilon^{10}}\,.
\nonumber
\end{eqnarray}
This is equivalent to
\begin{equation}
    u_t=-\epsilon u_x+u_{xx}
    +\frac{h^2}{12}(\epsilon-\partial_x)^2u_{xx}
    +\Ord{h^4}\,,
    \label{eq:cons}
\end{equation}
and so is indeed consistent to $\Ord{h^2}$ as $h\to 0$\,, independent
of $\epsilon$, with the advection-diffusion \pde~(\ref{eq:advec}).
The coefficient of the enhanced diffusion is the familiar value~$1$
when $\epsilon h$ is small.
However, when the advection speed is large enough (taking $\epsilon>0$
hereafter in this section for simplicity), the diffusion coefficient
$\nu_1$ increases to aid stabilisation.
The series for $\nu_1$ in~(\ref{eq:hiad}) is identical to that for
$({\epsilon h}/{2})\coth\left(\epsilon h/2\right)$, and numerical
summation of the series using the Shanks transform
\cite[\S8.1]{Bender81} suggests strongly that $\nu_1\sim \epsilon h/2$
as $\epsilon h\to \infty$ and is indeed within a few percent of this
value for $\epsilon h>4$\,, see Figure~\ref{fig:egcof}.
Thus for large advection speed $\epsilon$ on a finite width grid, the
centre manifold analysis promotes the model\footnote{If the advection
velocity is negative $\epsilon<0$, then various signs change and the
large $\epsilon h$ model remains an \emph{upwind} model, but is
consequently written in terms of forward differences.
This also occurs for the later model~(\ref{eq:back}) accurate to
errors $\Ord{\gamma^3}$\,.} (written in terms of the backward
difference operator $\nabla$)
\begin{equation}
    \dot u_j\approx -\frac{\epsilon}{h}\nabla u_j
    =-\frac{\epsilon}{h}\left( u_j-u_{j-1} \right)\,.
    \label{eq:fastad}
\end{equation}
This is not, and need not be, consistent with the
\pde~(\ref{eq:advec}) as $h\to 0$ because it applies for finite
$\epsilon h$.
That it should be relevant to~(\ref{eq:advec}) comes from the centre
manifold expansion in $\gamma$ albeit evaluated at $\gamma=1$ (via the
``holistic'' arrows in Figure~\ref{fig:conc}).
Exact solutions of~(\ref{eq:fastad}) are readily obtained.
For example, consider a point release from $j=0$ at $t=0$: $u_j(0)=1$
if $j=0$ and is $0$ otherwise.
Then the moment generating function
$G(z,t)=\sum_{j=0}^{\infty}z^ju_j(t)$ is easily seen to be that for a
Poisson probability distribution with parameter $\epsilon t/h$, namely
$G(z,t)=\exp[(z-1)\epsilon t/h]$.
Hence the mean location and variance of $u_j$ are
\begin{equation}
    \mu_j=\sigma_j^2=\frac{\epsilon t}{h}
    \quad\Rightarrow\quad
    \mu_x=\epsilon t
    \quad\mbox{and}\quad
    \sigma_x^2=\epsilon h t\,.
    \label{eq:meanvar}
\end{equation}
Thus for $\epsilon h$ \emph{not small}: this model has precisely the
correct advection speed~$\epsilon$; and although the variance is
quantitatively wrong, $\epsilon ht$ instead of~$2t$, at least it is
qualitatively correct for finite $\epsilon h$.
The centre manifold model~(\ref{eq:hiad}) is consistent for small $h$
and has the virtue of being always stable and will always maintain
non-negativity of solutions no matter how large the advection speed
$\epsilon$.

Second, explore the holistic finite difference model with
second-order interactions between adjacent elements, that is, the case
$\ell=3$ for which errors are $\Ord{\gamma^3}$, but to high order in the
advection $\epsilon$.
Computer algebra derives the model
\begin{eqnarray}
    \dot u_j & = &  -\frac{\gamma\epsilon}{h}\mu\delta u_j
    \nonumber\\&&{}
    +\frac{\gamma}{h^2}\left[ \delta^2
    +\frac{\epsilon^2h^2}{12}\delta^2
    -\frac{\epsilon^2h^2}{720}\delta^2
    +\frac{\epsilon^6h^6}{30240}\delta^2
    -\frac{\epsilon^8h^8}{1209600}\delta^2 \right] u_j
    \nonumber  \\
     &  & {}+\frac{\gamma^2}{h^2}\left[ - \frac{1}{12} \delta^4
     + \epsilon^2 h^2 \left( - \frac{1}{12} \delta^2
     - \frac{1}{30} \delta^4\right)
     + \epsilon^4 h^4 \left(\frac{1}{720} \delta^2
     + \frac{1}{5040} \delta^4\right)
    \right.\nonumber\\&&\left.\quad{}
     + \epsilon^6 h^6 \left( - \frac{1}{30240} \delta^2
     + \frac{1}{151200} \delta^4\right)
    \right.\nonumber\\&&\left.\quad{}
     + \epsilon^8 h^8 \left(\frac{1}{1209600} \delta^2
     - \frac{1}{1900800} \delta^4\right)
     \right]u_j
    \nonumber  \\
     &  & {}+\frac{\gamma^2\epsilon}{h}\left[ \frac{1}{6}
     + \frac{\epsilon^2 h^2}{90}
     - \frac{\epsilon^4 h^4}{2520}
     + \frac{\epsilon^6 h^6}{75600}
     - \frac{\epsilon^8 h^8}{2395008} \right]\mu\delta^3u_j
     \nonumber\\
     & & {}+\Ord{\gamma^3,\epsilon^{10}}\,.
    \label{eq:gamm}
\end{eqnarray}
Observe the marvellous feature that when we evaluate this at
$\gamma=1$ all the terms in $\delta^2$ associated with high orders of
$\epsilon h$ neatly cancel.
The model for the advection-diffusion thus reduces to
\begin{eqnarray}
    \dot u_j & = &
    -\frac{\epsilon}{h}\left(\mu\delta -\kappa_2{\mu\delta^3}\right) u_j
    +\frac{1}{h^2}\left({\delta^2} -\nu_2{\delta^4}\right) u_j\,,
    \label{eq:hiiad}  \\
     \mbox{where}\quad \nu_2 & = &
      \frac{1}{12}
     + \frac{\epsilon^2 h^2}{30}
     - \frac{\epsilon^4 h^4}{5040}
     - \frac{\epsilon^6 h^6}{151200}
     + \frac{\epsilon^8 h^8}{1900800}
     +\Ord{\epsilon^{10}}\,,
    \nonumber  \\
     \mbox{and}\quad \kappa_2 & = &
     \frac{1}{6}
     + \frac{\epsilon^2 h^2}{90}
     - \frac{\epsilon^4 h^4}{2520}
     + \frac{\epsilon^6 h^6}{75600}
     - \frac{\epsilon^8 h^8}{2395008}
     +\Ord{\epsilon^{10}}\,.
     \nonumber
\end{eqnarray}
See that in this model, as $h\to 0$, the hyperdiffusion coefficient
$\nu_2\sim 1/12$ and the dispersion coefficient $\kappa_2\sim 1/6$ to
give the classic second-order in $h$ corrections to the central
difference approximations of the first two derivatives.
Indeed the model~(\ref{eq:hiiad}) is equivalent to
\begin{equation}
    u_t=-\epsilon u_x+u_{xx}
    +\frac{h^4}{90}(\epsilon-\partial_x)^3u_{xxx}
    +\Ord{h^6}\,,
    \label{eq:conss}
\end{equation}
and so is consistent to $\Ord{h^4}$ as $h\to 0$\,, independent of
$\epsilon$, with the advection-diffusion \pde~(\ref{eq:advec}).

For large advection speed or grid size, large $\epsilon h$, the
model~(\ref{eq:hiiad}) is astonishingly good.
Using the large $\epsilon h$ approximations indicated by the Shanks
transforms plotted in Figure~\ref{fig:egcof} for $\nu_2$ and
$\kappa_2$, the model~(\ref{eq:hiiad}) reduces to simply
\begin{eqnarray}
    \dot u_j
    &=&-\frac{\epsilon}{h}\left( \nabla+\frac{1}{2}\nabla^2 \right)u_j
    +\frac{1}{h^2}\nabla^2 u_j
    \nonumber\\
    &=&-\frac{\epsilon}{2h}\left( u_{j-2}-4u_{j-1}+3u_j \right)
    +\frac{1}{h^2}\left( u_{j-2}-2u_{j-1}+u_j \right)\,.
    \label{eq:back}
\end{eqnarray}
This large $\epsilon h$ model uses only backward differences to
incorporate second-order upwind estimates of the derivative,
$\nabla+\frac{1}{2}\nabla^2$, and the second derivative, $\nabla^2$.
To show its good properties,\footnote{The upwind difference
model~(\ref{eq:back}) is only stable for $\epsilon h>2/3$\,.
However, from Figure~\ref{fig:egcof} the approximation~(\ref{eq:back})
is only applicable to~(\ref{eq:hiiad}) for $\epsilon h$ greater than
about $4$; thus its instability for very much smaller $\epsilon h$ is
irrelevant.} consider again a point release from $j=0$ at time $t=0$.
The moment generating function $G(z,t)=\sum_{j=0}^{\infty} z^ju_j(t)$
for the evolution governed by~(\ref{eq:back}) is readily shown to be
\begin{equation}
    G(z,t)=\exp\left[ -\frac{\epsilon t}{2h}(z-1)(z-3)
    +\frac{t}{h^2}(z-1)^2 \right]\,.
    \label{eq:mgf}
\end{equation}
Then since
\begin{displaymath}
    \left.\D zG\right|_{z=1}=\frac{\epsilon t}{h}
    \quad\mbox{and}\quad
    \left.\DD zG\right|_{z=1} =\left(\frac{\epsilon t}{h}\right)^2
    -\frac{\epsilon t}{h} +\frac{2t}{h^2}\,,
\end{displaymath}
we determine the mean position and variance of the spread in $u_j$ to
be
\begin{equation}
    \mu_j=\frac{\epsilon t}{h}\mbox{ and }\sigma_j^2=\frac{2t}{h^2}
    \quad\Rightarrow\quad
    \mu_x=\epsilon t\mbox{ and }\sigma_x^2=2t\,.
    \label{eq:menvar}
\end{equation}
This predicted mean and variance following a point release are
\emph{exactly} correct \emph{for all} time for the advection-diffusion
\pde~(\ref{eq:advec}).  This specific result applies to all finite
advection speeds $\epsilon$ and all finite grid spacings $h$
whenever $\epsilon h$ is large enough.

It seems that creating finite differences which, as shown in
Figure~\ref{fig:conc}, are both consistent for small grid spacing~$h$
and also holistically derived via centre manifold theory thus can lead
to models which are remarkably accurate and stable over a wide range
of parameters.
The hierarchy of refinements in the coefficient, $\epsilon$, of odd
derivatives, here just $\partial_x$, seem to be useful for robust
performance for finite~$h$.

\section{Models of the even terms are consistent}
\label{Sline}

In this section I prove that the proposed centre manifold approach
consistently models all the even derivatives in $\cA$.
That is, the equivalent \pde{} of the finite difference model on the
centre manifold~$\cM$ matches $u_t=\cA u$ to some order in grid
spacing~$h$; the order of error is controlled by the order of
truncation, $\ell$, in the coupling coefficient $\gamma$.
The veracity of the following theorems are supported by the
results of suitable variants of the computer algebra program of
Appendix~\ref{Scompalg}.

Remarkably, the polynomials found here to describe the structure
within each element are \emph{independent} of the linear operator
$\cA$!

\begin{theorem}\label{thm:cm}
The centre manifold model~(\ref{eq:mod}) constructed with the
amplitude condition~(\ref{eq:amp}), the internal boundary
condition~(\ref{eq:ibcv}) and to errors~$\Ord{\gamma^\ell}$, forms a
sem-discrete finite difference approximation to the \pde{} $u_t=\cA u$
consistent to~$\Ord{\partial_x^{2\ell}}$, where $\cA$ is any even
operator.
\end{theorem}

\begin{proof}
I construct the proof in stages beginning with the end result and
finding successive sufficient conditions for the preceding steps.
Observe that I actually prove a slightly stronger result: by allowing
the even operator $\cA$ to contain a constant term $a_0$,
see~(\ref{eq:cl}), I prove the consistency of an invariant manifold
model based upon the mode with eigenvalue $\lambda_0=a_0$.
When $a_0\geq 0$ this forms a centre-unstable manifold model for which
a relevance theorem also ensures asymptotic completeness.
\begin{itemize}
    \item To solve $u_t=\cA u$ to errors $\Ord{\gamma^\ell}$, expand
    the centre manifold model~(\ref{eq:mod}) to~$\Ord{\gamma^\ell}$ as
        \begin{equation}
        u=v(\vec u,x)=\sum_{k=0}^{\ell-1}\gamma^kv^k(\vec u,x)\,,
	\quad\mbox{and}\quad
	\dot u_j=g(\vec u)=\sum_{k=0}^{\ell-1}\gamma^kg^ku_j\,,
        \label{eq:lexp}
    \end{equation}
    where $g^k$ are some difference operators, and where $v$ and $g$
    implicitly refer to the $j$th element.  Note that superscripts
    to~$\gamma$ denote exponentiation whereas those on~$v$ and~$g$
    denote the index of coefficients in the asymptotic expansion: I do
    this because~$v$ and~$g$ have an implicit subscript~$j$ denoting
    the grid element under consideration.  Then substitution into the
    \pde{} and extracting powers of $\gamma$ shows that we require
        \begin{equation}
        \cA v^k=g^0v^{k}+g^1v^{k-1}+\cdots+g^{k-1}v^1+g^kv^0
	\quad\mbox{for $0\leq k<\ell$\,.}
        \label{eq:ordk}
    \end{equation}
    Similarly, substitution into the amplitude condition~(\ref{eq:amp})
    requires
        \begin{equation}
        v^0(\vec u,x_j)=u_j\,,
	\quad\mbox{and}\quad
	v^k(\vec u,x_j)=0\quad\mbox{for $1\leq k<\ell$\,.}
        \label{eq:ampk}
    \end{equation}
    Whereas substitution into the \ibc~(\ref{eq:jbc}) requires,
    evaluating the left-hand sides at $x_j$,
        \begin{eqnarray}
        \mu_x\delta_x v^k & = & \left\{ \begin{array}{ll}
            \mu\delta u_j\,, & k=1\,,  \\
            0\,, & k\neq 1\,,
        \end{array}\right.
        \label{eq:ibcak}  \\
	\mbox{and}\quad
        \delta_x^2 v^k&=&\left\{
                \begin{array}{ll}
            \delta^2 u_j\,, & k=1\,,  \\
            0\,, & k\neq 1\,.
        \end{array}\right.
        \label{eq:ibcbk}
    \end{eqnarray}
    Equations~(\ref{eq:ordk}--\ref{eq:ibcbk}) form a well-posed system
    of equations for $v^k$ and~$g^k$.  In many applications of centre
    manifold theory, because $\cA-g^0$ is singular, we often solve
    each level in the hierarchy of equations in two steps: the first
    is to find $g^k$ by ensuring that the right-hand side
    of~$\ref{eq:ordk}$ is in the range of $\cA-g^0$, this is the
    so-called ``solvability condition''; the second step is to find
    $v_k$.  However, here we proceed to construct straightforwardly
    the solution of the entire set of equations in general.

\item We show the hierarchy of differential equations~(\ref{eq:ordk})
are satisfied by functions~$v^k$ and give \emph{consistent} finite
difference operators~$g^k$, if $v^k$ satisfy the recursive difference
equation
        \begin{equation}
        \delta_x^2v^k=\delta^2v^{k-1}\quad\mbox{for all $x$,}
	\quad\mbox{and}\quad
	v^0=\mbox{constant}\,.
        \label{eq:diffk}
    \end{equation}
        \begin{itemize}
	\item By the following induction argument~(\ref{eq:diffk})
	implies that
                \begin{equation}
            \delta_x^{2m}v^k=\left\{
                        \begin{array}{ll}
                \delta^{2m}v^{k-m}\,,&\mbox{for }m=0,\ldots,k\,,  \\
                0\,, & m=k+1,k+2,\ldots\,.
            \end{array}\right.
            \label{eq:difffk}
        \end{equation}
	Now,
	    \begin{eqnarray*}
	    \delta_x^{2m}v^k & = & \delta_x^{2m-2}\delta_x^2v^k  \\
	     & = & \delta_x^{2m-2}\delta^2v^{k-1}
	     \quad\mbox{by (\ref{eq:diffk})}\\
	     & = & \delta^2\delta_x^{2(m-1)}v^{k-1}
	     \quad\mbox{as $\delta$ and $\delta_x$ commute}\\
	     & = & \delta^2\delta^{2(m-1)}v^{k-m}
	     \quad\mbox{p.v.~(\ref{eq:difffk}) holds for $m-1$}\\
	     & = & \delta^{2m}v^{k-m}\,.
	\end{eqnarray*}
	Then since~(\ref{eq:difffk}) is trivially true for $m=0$, it
	follows by induction that~(\ref{eq:difffk}) holds for all $m$
	provided $m\leq k$.  Further, since $v^0$ is constant
	by~(\ref{eq:diffk}), $\delta_x^{2k}v^k$ is constant, so higher
	order differences ($m>k$) are all zero.

	\item By an ``even'' operator $\cA$ I mean one which only
	involves even order derivatives in $x$.  Hence write $\cA$
	formally as an infinite sum of even powers of the central
	difference operator
                \begin{equation}
            \cA=\sum_{m=0}^{\infty}a_m\delta_x^{2m}\,,
            \label{eq:cl}
        \end{equation}
	for some coefficients $a_m$; for example, for the diffusion
	operator in~(\ref{eq:advec}), from~\cite[p65,e.g.]{npl61},
	    \begin{displaymath}
	    \DD x{}=\frac{4}{h^2}\arcsinh^2(\half\delta_x)
	    =\frac{1}{h^2}\delta_x^2 -\frac{1}{12h^2}\delta_x^4
	    +\frac{1}{90h^2}\delta_x^6-\cdots\,;
	\end{displaymath}
	more generally, $\cA$ could be any symmetric convolution
	operator for which the infinite sum~(\ref{eq:cl}) forms a
	reasonable representation.  Then~(\ref{eq:difffk})
	ensures~(\ref{eq:ordk}) since
	    \begin{eqnarray*}
	    \cA v^k & = & \sum_{m=0}^{\infty}a_m\delta_x^{2m}v^k  \\
	     & = & \sum_{m=0}^k a_m\delta^{2m}v^{k-m}
	     \quad\mbox{by (\ref{eq:difffk})}\\
	     & = & g^0v^{k}+g^1v^{k-1}+\cdots+g^{k-1}v^1+g^kv^0\,,
	\end{eqnarray*}
	provided $g^k=a_k\delta^{2k}$ which are \emph{precisely the
	operators required to make the model $\dot u_j=g(\vec u)$ of
	$u_t=\cA u$ consistent to~$\Ord{\partial_x^{2\ell}}$, when
	truncated as in~(\ref{eq:lexp}) to
	errors~$\Ord{\gamma^\ell}$.}
    \end{itemize}

    \item By simple substitution, a sequence of functions $v^k$
    satisfying the recurrence~(\ref{eq:diffk}), amplitude
    conditions~(\ref{eq:ampk}) and the internal boundary
    conditions~(\ref{eq:ibcak}--\ref{eq:ibcbk}) are
        \begin{equation}
        v^0=u_j\,,\quad
	v^k=p_k(\xi)\mu\delta^{2k-1}u_j +q_k(\xi)\delta^{2k}u_j
	\,,\quad\mbox{for }k\geq1\,,
        \label{eq:vsol}
    \end{equation}
    where, as always, $\xi=(x-x_j)/h$ is a grid scaled coordinate and
    provided that for $k\geq 1$
        \begin{equation}
        \delta_x^2p_k=p_{k-1}\,, \quad
	p_k(\pm k)=\pm1\,,\quad
	\mbox{and}\quad
	p_k(\xi)=0\mbox{ for }\xi=0,\pm 1\,,
        \label{eq:p}
    \end{equation}
    and similarly
            \begin{equation}
        \delta_x^2q_k= q_{k-1}\,, \quad
	q_k(\pm k)=+\half\,,\quad
	\mbox{and}\quad
	q_k(\xi)=0\mbox{ for }\xi=0,\pm 1\,,
        \label{eq:q}
    \end{equation}
    after defining $q_0(x)=1$ and $p_0(x)=0$.

    \item Analysing the difference tables for $p_k$ and $q_k$ and
    straightforward induction using~(\ref{eq:p}--\ref{eq:q}) proves
    that $p_k(\xi)=q_k(\xi)=0$ for integer $\xi\in[-k+1, k-1]$.  Then
    the following $p_k$ and $q_k$ are the unique polynomials, of
    degree~$2k-1$ and~$2k$ respectively, also satisfying $p_k(\pm
    k)=\pm1$ and $q_k(\pm k)=+\half$:
        \begin{equation}
        p_k(\xi)=\frac{1}{(2k-1)!}\prod_{m=-k+1}^{k-1}(\xi-m)\,,
	\quad\mbox{and}\quad
	q_k(\xi)=\frac{\xi}{2k}p_k(\xi)\,,
        \label{eq:pksol}
    \end{equation}
    as plotted in Figure~\ref{fig:pq}.  For example, $p_1(\xi)=\xi$
    and $q_1(\xi)=\half \xi^2$.
\begin{figure}[tbp]
    \centering
    \includegraphics{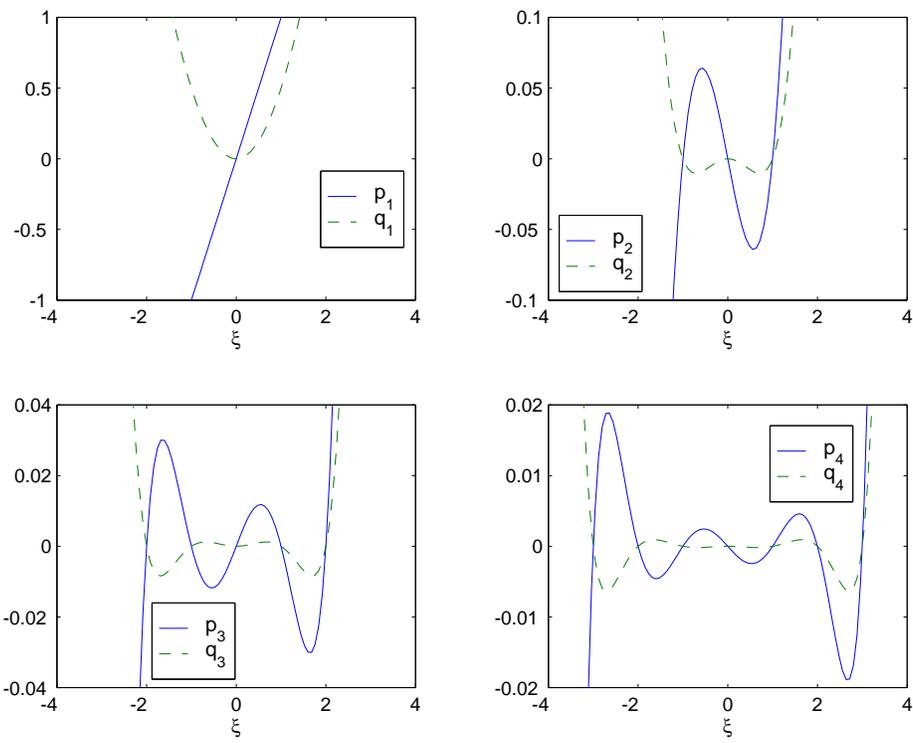}
    \caption{graphs of the polynomials~(\ref{eq:pksol}) forming a
    basis for the fields of the approximations to the centre
    manifold.}
    \label{fig:pq}
\end{figure}

    These polynomial $p_k$ and $q_k$ also need to satisfy the
    recurrences in~(\ref{eq:p}--\ref{eq:q}) \emph{pointwise} in $\xi$.
    This is trivially true for $k=1$.  Now, $\delta_x^2p_{k+1}$ is
    from~(\ref{eq:pksol}) a polynomial of degree~$2k-1$, from its
    difference table has the same zeros as~$p_k$, it is $\pm1$ at $\pm
    k$, and so must be $p_k(\xi)$ for all $\xi$.  Similarly for
    $\delta_x^2q_{k+1}$.\footnote{One easily imagines other
    functions~$p_k$ and~$q_k$ that have all the requisite properties
    to ensure a consistent finite difference approximation $g(\vec
    u)$.  For example, $\tilde p_k(\xi)=p_k(\xi)+a\sin(2n\pi \xi)$\,.
    However, as is often the case, if the method of construction of
    the centre manifold makes $v^k$ a polynomial of degree~$2k$, then
    the solution for $v^k$ is the one given in~(\ref{eq:vsol}) in
    terms of $p_k$ and $q_k$.  }
\end{itemize}
\end{proof}

\begin{corollary}\label{cor:obv}
It follows immediately from the theorem that if the highest derivative
in the even operator $\cA$ is of order~$n$, then the finite difference
model~(\ref{eq:lexp}) is consistent with $u_t=\cA u$ to
$\Ord{h^{2\ell-n}}$ as $h\to 0$.
\end{corollary}

\section{Odd perturbations are also consistent}
\label{Slino}

The results of the previous section on the modelling of even operators
$\cA$ are extremely satisfactory.  The modelling of odd operators is
not quite so neat.  In the general linear \pde~(\ref{eq:pde}) I
introduced the odd terms, $\cB$, with a multiplying $\epsilon$.  The
reason is, as seen in \S\ref{Sadvec}, that such odd terms generate
extra terms in the finite difference model which are nonlinear in the
coefficients of the odd derivatives, that is $\Ord{\epsilon^2}$.  As
ellaborated in \S\ref{Sadvec} for the advection-diffusion equation,
these higher-order contributions seem to reflect the changes needed
for stable discretisations of equations with large amounts of
advection, $u_x$, or dispersion, $u_{xxx}$.  The extra complications
of these nontrivial effects of odd derivatives appear necessary.
However, here we restrict attention to proving consistency to an error
quadratic in the odd coefficients and leave to further research the
investigation of higher-order consistency.

\begin{theorem}\label{thm:odd}
The centre manifold model~(\ref{eq:mod}) constructed with the
amplitude condition~(\ref{eq:amp}), the internal boundary
condition~(\ref{eq:ibcv}) and to
errors~$\Ord{\gamma^\ell,\epsilon^2}$, forms a finite difference
approximation to the \pde{} $u_t=\cA u+\epsilon\cB u$ consistent
to~$\Ord{\partial_x^{2\ell} +\epsilon\partial_x^{2\ell-1},
\epsilon^2}$, where $\cA$ is any even operator and $\cB$ is any odd
operator.
\end{theorem}

\begin{proof}
As in (\ref{eq:lexp}), we expand the the centre manifold
ansatz~(\ref{eq:mod}) to errors $\Ord{\gamma^\ell,\epsilon^2}$:
    \begin{equation}
        u       =\sum_{k=0}^{\ell-1}\gamma^kv^k
	+\epsilon\sum_{k=0}^{\ell-1}\gamma^kw^k\,,
	\quad\mbox{and}\quad
	\dot u_j=\sum_{k=0}^{\ell-1}\gamma^kg^ku_j
	+\epsilon\sum_{k=0}^{\ell-1}\gamma^kf^ku_j\,,
        \label{eq:lexpo}
    \end{equation}
where $f^k$ and $g^k$ are some difference operators, and $v^k$ and
$w^k$ are functions defined about the $j$th gridpoint, and they all
implicitly refer to the $j$th element (as before the superscript to
$v$, $w$, $g$ and~$h$ denotes an index in the series, whereas the
superscript to~$\gamma$ denotes exponentiation).  After substitution
into the \pde{}, terms in $\epsilon^0$ determines $v^k$ and $g^k$ as
in the previous section.  Upon extracting from the \pde{} the
coefficients of terms linear in $\epsilon$ and of various powers of
$\gamma$ requires us to solve
     \begin{equation}
        \cA w^k+\cB v^k=\sum_{r=0}^{k}\left(f^{k-r}v^{r}
	+g^{k-r}w^{r}\right)\,,
	\quad\mbox{for $0\leq k<\ell$\,.}
        \label{eq:ordko}
    \end{equation}
Substitution of the expansion~(\ref{eq:lexpo}) into the amplitude
condition~(\ref{eq:amp}) and the \ibc's~(\ref{eq:jbc})
and equating coefficients of $\epsilon\gamma^k$ leads to these
internal boundary conditions for the $w^k(x)$:
        \begin{equation}
        w^k=0\quad\mbox{for $x=x_j,x\jpm$}\,.
        \label{eq:wibc}
    \end{equation}
Since $\cB$ is an odd operator, we formally write it as the
following infinite sum of odd powers of centred difference
operators
\begin{equation}
        \cB=\sum_{m=1}^{\infty} b_m\mu_x\delta_x^{2m-1}\,,
        \label{eq:oddb}
    \end{equation}
for some coefficients $b_m$; for example (from~\cite[p65,e.g.]{npl61})
        \begin{displaymath}
        \D x{}=\frac{2\mu_x}{h}\arcsinh(\half\delta_x)
	=\frac{1}{h}\mu_x\delta_x -\frac{1}{6h}\mu_x\delta_x^3
	+\frac{1}{30h}\mu_x\delta_x^5 -\cdots\,;
\end{displaymath}
and more generally, $\cB$ could be any antisymmetric convolution
operator for which the infinite sum~(\ref{eq:oddb}) is reasonable.
I prove that there exists solutions $w^k(x)$, odd functions of $x$
(about $x_j$), with
\begin{equation}
    f^k=b_k\mu\delta^{2k-1}\,,
    \label{eq:fk}
\end{equation}
so that the model~(\ref{eq:lexpo}) is consistent with the effect of
the odd derivatives in $\epsilon\cB$ to errors
$\Ord{\delta_x^{2\ell-1}} =\Ord{\partial_x^{2\ell-1}}$.
Since we already know $v^k$ and since $w^k$ appears to vary for
different problems, the operators $f^k$ are determined by the
solvability condition that all terms except $\cA w^k-g^0 w^k$
appearing in~(\ref{eq:ordko}) combine to be in the range of the
singular $\cA-g^0$.

\begin{itemize}
    \item First, prove that the even part of $\cB v^k$ cancels with
    the even part of $\sum_{r=0}^{k}f^{k-r}v^{r}$ and so is eliminated
    from~(\ref{eq:ordko}).  As a preliminary step consider,
    using~(\ref{eq:pksol}),
        \begin{eqnarray}
        \mu_x\delta_xp_k(\xi) & = & \frac{1}{(2k-1)!}
        \frac{1}{2}\left[ \prod_{m=-k+1}^{k-1}(\xi+1-m)
        -\prod_{m=-k+1}^{k-1}(\xi-1-m)\right]
	\nonumber\\
         & = & \frac{1}{(2k-1)!}\prod_{m=-k+2}^{k-2}(\xi-m)
	 \times\nonumber\\&&\quad
         \frac{1}{2}\left[ (\xi+k)(\xi+k-1)-(\xi-k)(\xi-k+1) \right]
	 \nonumber\\
         & = & \frac{1}{(2k-1)(2k-2)}p_{k-1}(\xi)\times(2k-1)\xi
	 \nonumber\\
         & = & q_{k-1}(\xi)\,.\label{eq:mdp}
    \end{eqnarray}
    Then from~(\ref{eq:vsol}) and since $p_k$ is odd and $q_k$ is
    even, see Figure~\ref{fig:pq}, $\cB q_k$ is odd and so the even
    part of $\cB v^k$ is
        \begin{eqnarray}
        Bp_k\mu\delta^{2k-1}u_j & = & \sum_{m=1}^{\infty}
        b_m\mu_x\delta_x^{2m-1} p_k \,\mu\delta^{2k-1}u_j
        \quad\mbox{by (\ref{eq:oddb})}\nonumber  \\
         & = & \sum_{m=1}^{k}
        b_m\mu_x\delta_x p_{k-m+1} \,\mu\delta^{2k-1}u_j
        \quad\mbox{by (\ref{eq:q}) inductively}\nonumber  \\
         & = & \sum_{m=1}^{k}
        b_m q_{k-m} \,\mu\delta^{2k-1}u_j
        \quad\mbox{by (\ref{eq:mdp}).}
        \label{eq:bpk}
    \end{eqnarray}
    Whereas on the right-hand side of~(\ref{eq:ordko}) the even part
    of $\sum_{r=0}^{k}f^{k-r}v^{r}$ is
    \begin{eqnarray}
        \sum_{r=0}^{k}f^{k-r}q^{r}\delta^{2r}u_j & = &
	\sum_{r=0}^k q^{r}b_{k-r}\mu\delta^{2(k-r)-1} \delta^{2r}u_j
	\quad\mbox{by (\ref{eq:fk})} \nonumber\\
         & = & \sum_{r=0}^k q^{r}b_{k-r}\,\mu\delta^{2k-1} u_j\,,
        \label{eq:mdpr}
    \end{eqnarray}
    which exactly cancels with~(\ref{eq:mdp}) from the even part of
    $\cB v^k$ on the left-hand side of~(\ref{eq:ordko}).

    \item Second, since the even components of~(\ref{eq:ordko}) that
    involve the various $v^k$ cancel, and since $\cA$ and $g^k$ are
    even, then a particular solution $w^k(x)$ of~(\ref{eq:ordko}) may
    be found that is odd.  The conditions~(\ref{eq:wibc}) then force
    these odd $w^k$ to be the unique solutions in the space of finite
    degree polynomials.

    \item Third, consider evaluating the hierarchy of
    equations~(\ref{eq:ordko}) at the centre grid point of the
    element, $x=x_j$ or equivalently $\xi=0$, and simplify the various
    contributions.
\begin{itemize}
    \item  Now $\cA$ is an even operator and $w^k$ an odd function,
    so $\cA w^k$ is an odd function, and thus $\cA w^k(x_j)=0$.

    \item Since $g^{k-r}$ is a difference operator in~$j$ it does not
    affect the amplitude condition that $w^r(x_j)=0$, and thus
    $g^{k-r}w^r(x_j)=0$.

    \item I have already shown that the even parts of $\cB v^k$ and
    $\sum_{r=0}^{k}f^{k-r}v^r$ agree pointwise, so they certainly do
    at $x_j$, at which the odd parts must also trivially vanish
    together.
\end{itemize}
Thus the choice~(\ref{eq:fk}) is the unique one to satisfy this
solvability condition for~(\ref{eq:ordko}).

\end{itemize}
Since the expansion~(\ref{eq:lexpo}) then contains the exact
differences up to $b_{\ell-1}\mu\delta^{2\ell-3}$ and
$a_{\ell-1}\delta^{2\ell-2}$, the errors in the finite truncation of
the finite difference model will be $\Ord{\partial_x^{2\ell}}$ from
the even terms and $\Ord{\epsilon\partial_x^{2\ell-1},\epsilon^2}$
from the odd.

\end{proof}

\begin{corollary}\label{cor:obvo}
    It follows immediately from the theorem that if the highest
    derivative in the operator $\cA+\epsilon\cB$ is of order~$n$, then
    the finite difference model~(\ref{eq:lexp}) is consistent with
    $u_t=\cA u+\epsilon\cB u$ to $\Ord{h^{2\ell-1-n}}$ as $h\to0$ to
    an error $\Ord{\epsilon^2}$.
\end{corollary}

\section{Concluding remarks}

We have seen that the artificial internal boundary
conditions~(\ref{eq:ibcv}) together with the application of centre
manifold theory generate finite difference models that have remarkably
good properties, at least for linear systems.  These are explcitly
shown for the example of the advection-diffusion equation
(\S\ref{Sadvec}) where we saw not only consistency for small $h$, but
also an appropriate upwind model for large $\epsilon h$.  Although the
Theorem~\ref{thm:odd} only establishes consistency with the odd terms
to $\Ord{\epsilon^2}$ the advection-diffusion example of
\S\ref{Sadvec} shows that higher-order consistency is possible.
Further research is needed on the characteristics of higher orders in
the odd derivatives.

Also, further research, such as that for Burgers' equation
in~\cite{Roberts98a}, will explore the performance of this holistic
approach to discretisations of nonlinear systems in various spatial
dimensions.  Since centre manifold theory is designed to analyse
nonlinear systems I expect reliable models to be derived.

Throughout the analysis in this paper I have parametrised the centre
manifold model in terms of the field at each of the grid points,
$u_j=u(x_j,t)$.  This was done for simplicity.  Other choices are
possible for the parameters of the finite difference model, for
example we could choose to use the element average
\begin{equation}
    u_j(t)=\frac{1}{h}\int_{x_j-h/2}^{x_j+h/2} u(x,t)\,dx\,.
    \label{eq:ampav}
\end{equation}
This choice would be appropriate to easily establish the conservation
of total $u$, or not as the case may be.  Computational experiments
show that either of these choices of amplitude produce equivalent
results for linear systems.

Centre manifold theory is routinely applied to autonomous dynamical
systems.  However, the geometric viewpoint it establishes leads to a
rational treatment of the projection of initial conditions onto the
finite dimensional model and of the projection of a perturbing
forcing \cite{Roberts97b, Cox93b, Roberts89b}.

\paragraph{Acknowledgement:} this research was supported by a grant
from the Australian Research Council.

\appendix
\section{Computer algebra for perturbed diffusion}
\label{Scompalg}

This computer algebra program is included to replace the recording of
tedious details of elementary algebra involved in solving the
advection-diffusion \pde~(\ref{eq:advec}).

For this problem the iterative algorithm of~\cite{Roberts96a} is
implemented in \textsc{reduce}\footnote{At the time of writing,
information about \texttt{reduce} was available from Anthony C.~Hearn,
RAND, Santa Monica, CA~90407-2138, USA.
\url{mailto:reduce@rand.org}} Although there are many details in the
program, the correctness of the results are \emph{only determined} by
driving to zero (line~40) the residuals of the governing \pde,
evaluated on line~31, and the residuals of the \ibc's, evaluated on
lines~32--35, to the error specified on line~29.
The other details in the program only affect the rate of convergence
to the ultimate answer.

Lines~44--46 then rewrite the model in terms of the central difference
operator~$\delta$ and mean~$\mu$.  Lastly, lines~48--53 derive the
equivalent differential equation for the finite difference model.

Other problems in the same class can be and have been examined by
amending line~31, the evaluation of the residual of the governing
\pde, with other differential operators such as
\begin{quote}
\begin{verbatim}
+eps*(a0*v+a2*df(v,x,4)+b2*df(v,x,3))
\end{verbatim}
\end{quote}
The results of the analysis with such terms confirm the theorems in
Sections~\ref{Sline} and~\ref{Slino}.

\verbatimlisting{linear.red}

\bibliographystyle{plain}
\bibliography{ajr,bib,new}

\end{document}